\documentclass[12pt]{article}
\usepackage[utf8]{inputenc}
\usepackage{cmap}
\usepackage{amsmath}
\usepackage{amssymb}
\usepackage{amscd}
\usepackage[english]{babel}
\usepackage{graphicx}
\usepackage{url}
\begin{document}
\author{V.~K.~Beloshapka\thanks{Faculty of Mechanics and Mathematics,
Lomonosov Moscow State University, Vorobyovy Gory, 119991 Moscow, Russia;
vkb@strogino.ru}}

\newcommand{\CC}{\mathbb C}
\newcommand{\RR}{\mathbb R}
\newcommand{\Ree}{\operatorname{Re}}
\newcommand{\Imm}{\operatorname{Im}}
\newcommand{\aut}{\mathfrak{aut}}
\newcommand{\Span}{\operatorname{span}}
\newcommand{\Dz}[1]{D_{#1}}
\newcommand{\Ew}[1]{E_{#1}}

\date{July 21, 2026}

\title{\bf On Exceptional $CR$-Quadrics: Further Developments}

\maketitle

\begin{abstract}
Exceptional $CR$-quadrics are studied. An example of an exceptional
$(4,4)$-quadric is constructed; it realizes the minimum exceptional type
with respect to both $n$ and $k$. Its graded Lie algebra is described.
The available information on exceptional types is summarized, and the
lattice of $CR$-types is decomposed into the disjoint union of three sets:
$A$, the types for which exceptional quadrics are impossible; $B$, the
types for which examples of exceptional quadrics are known; and $C$, the
types whose status is currently unknown (neither an example nor a
nonexistence proof is known). Several questions are posed.
\end{abstract}

{\bf 1. Introduction}

\vspace{3ex}

This paper continues \cite{VB21}, where the study of exceptional
$CR$-quadrics was initiated. Paper \cite{VB21} was motivated by Meylan's
work \cite{FM}, which gave the first example of a quadric for which the
positive part of the graded Lie algebra of holomorphic automorphisms has
length greater than two, and by its continuation, the paper of Meylan and
Gregorovi\v{c} \cite{GM}, where a series of such examples was constructed.
Exceptional quadrics were also studied by Tumanov in \cite{T}. We continue
the study of this phenomenon, using results from \cite{FM}, \cite{GM}, and
\cite{T}. For general background on rigidity and finite jet determination
in $CR$-geometry, see \cite{CM,BER,Zaitsev}; for graded Lie algebras and
prolongation, see \cite{Tan,Yamaguchi,CapSlovak}.

\vspace{3ex}

We recall the basic definitions \cite{VB91,VB02}. Fix positive integers
$n$ and $k$ such that $1\leq n$ and $1\leq k\leq n^2$. Let
$(H_1,\dots,H_k)$ be a collection of Hermitian matrices of order $n$.
The collection is called {\it nondegenerate} if the matrices are linearly
independent and have no common kernel, that is,
$H_1z=\dots=H_kz=0$ implies $z=0$. Such a collection determines a
$CR$-quadric of type $(n,k)$, where $n$ is the dimension of the complex
tangent space and $k$ is the codimension. Let $z\in\mathbf C^n$,
$w\in\mathbf C^k$, and put
$<z,\bar z>=(H_1z\cdot\bar z,\dots,H_kz\cdot\bar z)$.
Then $Q$ is the submanifold of $\mathbf C^{n+k}$ defined by
$$ \mathrm{Im}\,w =<z,\bar{z}>.$$  
Let $\mathrm{Aut}\,Q$ be the local group of holomorphic automorphisms of
the germ of $Q$ at the origin, and let $\mathrm{aut}\,Q$ be its Lie
algebra. Nondegeneracy of the collection is the criterion for
finite-dimensionality of $\mathrm{aut}\,Q$. Introduce a grading on this
algebra by assigning weight $1$ to the variables $z$, weight $2$ to the
variables $w$, and the opposite weights to the corresponding
differentiations. Then
\begin{equation*}
 \mathrm{aut}\,Q=g_{-2}+g_{-1}+g_{0}+g_{1}+\dots+g_{l},   \quad  l < \infty,  
\end{equation*}
where $g_j$ is the component of weight $j$. The dimensions of $g_{-2}$
and $g_{-1}$ do not depend on $<z,\bar z>$ and are equal to $k$ and $2n$,
respectively. By contrast, the dimensions of the nonnegative components,
as well as the length $l$ of the algebra, depend on the Hermitian form.
Before \cite{FM}, the examples arising in this area had only the values
$0$, $1$, and $2$ for $l$. This holds, for example, when $k=1$ or $2$.
It was proved in \cite{VB21} that it also holds for $k=3$, and for
$k=4$, $n\leq3$.
A quadric is called {\it exceptional} if $l\geq3$. The first published
example of an exceptional quadric, given in \cite{FM}, had $CR$-type
$(4,5)$. An exceptional quadric of type $(6,4)$ and an infinite series of
examples of various types were constructed in \cite{GM}.

\vspace{3ex}

{\bf 2. An Exceptional Quadric of Type $(4,4)$}

\vspace{3ex}

It follows from the preceding discussion that codimension $4$ is the
smallest codimension in which exceptional quadrics may exist. It remained
open, however, whether $CR$-dimension $6$ was minimal when $k=4$; see
\cite{VB21}. The answer is negative. We now give an exceptional quadric of
type $(4,4)$. Neither parameter of this type can be decreased.

In $\mathbf C^8$ with coordinates
$(z_1,z_2,z_3,z_4,w_1,w_2,w_3,w_4)$, consider the nondegenerate quadric
$Q=Q_{44}$ of type $(4,4)$
\begin{equation*}
\begin{aligned}
\mathrm{Im} \, w_1
&=z_1\overline z_2+z_2\overline z_1,\\
\mathrm{Im} \,w_2
&=iz_1\overline z_2-iz_2\overline z_1,\\
\mathrm{Im} \,w_3
&=|z_2|^2,\\
\mathrm{Im} \,w_4
&=z_1\overline z_4+z_4\overline z_1
  +z_2\overline z_3+z_3\overline z_2.
\end{aligned}
\end{equation*}

\vspace{3ex}

This $(4,4)$-quadric is essentially a slight modification of Meylan's
exceptional $(4,5)$-quadric \cite{FM}. From the viewpoint of the
exceptionality criterion in \cite{VB21}, the decisive feature in the
construction of both examples is the existence of a linear syzygy. If
$(H_1,H_2,H_3,H_4)$ are the matrices of the coordinate forms, this syzygy
is the relation
\[i \, z_2 \, H_1\,z+  z_2 \, H_2\,z- 2\, i \, z_1 \, H_3\,z=0.\] 

The pencil corresponding to $Q$ is singular, that is,
\[P_H(t)= \mathrm{det}\,(t_1\,H_1+t_2\,H_2+t_3\,H_3+t_4\,H_4) \equiv 0.\]
Such quadrics are usually called {\it zero quadrics}. Note that Theorem~6
of \cite{VB21} immediately implies that a quadric of type $(n,n)$ which is
not a zero quadric is not exceptional.

\vspace{3ex}

Put
\[
\Dz{\alpha}:=\frac{\partial}{\partial z_\alpha},\qquad
\Ew{j}:=\frac{\partial}{\partial w_j} \qquad  1 \leq \alpha, \, j \leq 4.
\]
If $X$ is one of the holomorphic vector fields written below, then its
real part $2\,\mathrm{Re}\,X$ is an element of the real Lie algebra
$\mathrm{aut}(Q)$.

For brevity, put
\begin{eqnarray*}
 W_{+}:=w_1+iw_2,\qquad W_{-}:=w_1-iw_2,\\
 T=(-w_1+iw_2)z_2+2w_3z_1=2w_3z_1-W_{-} z_2.
\end{eqnarray*}
Since the negative-weight components have a common description for all
nondegenerate quadrics, only the nonnegative-weight components are listed
below.

\vspace{3ex}

{\bf Proposition 1.} (a) The graded Lie algebra $\mathrm{aut}\,Q$ has the
following structure:
$\mathrm{aut}\,Q=
g_{-2}\oplus g_{-1}\oplus
 g_0\oplus  g_1\oplus g_2
\oplus g_3\oplus  g_4$,
where
\begin{equation*}
\begin{array}{c|rrrrrrr}
j&-2&-1&0&1&2&3&4\\
\hline
\mathrm{dim}\, g_j&4&8&14&14&11&6&1
\end{array}  \;.
\end{equation*}
The coefficients of the vector fields are polynomials of degree at most
three.\\
(b) The local group $\mathrm{Aut}\,Q$ generated by $\mathrm{aut}\,Q$ is a
$58$-dimensional subgroup of the group of birational transformations of
$\mathbf C^8$, and the degrees of its elements do not exceed three. Its
closure in the Cremona group is an algebraic Lie subgroup.\\
{\it Proof.}
For a nondegenerate model quadric, the graded Lie algebra of its
infinitesimal $CR$-automorphisms is naturally identified with the full
Tanaka prolongation; see \cite{Tan,GM}. Computing the full Tanaka
prolongation of the symbol of $Q$, we obtain the following vector fields.

A basis of $g_0$:
\begin{align*}
X_{0,1}&=-z_1\Dz3+z_2\Dz4,\\
X_{0,2}&=iz_2\Dz3,\\
X_{0,3}&=iz_1\Dz4,\\
X_{0,4}&=iz_1\Dz3+iz_2\Dz4,\\
X_{0,5}&=i(z_1\Dz1+z_2\Dz2+z_3\Dz3+z_4\Dz4),\\
X_{0,6}&=z_2\Dz1-z_4\Dz3+2w_3\Ew1,\\
X_{0,7}&=iz_2\Dz2+iz_3\Dz3-w_2\Ew1+w_1\Ew2,\\
X_{0,8}&=z_1\Dz1-z_4\Dz4+w_1\Ew1+w_2\Ew2,\\
X_{0,9}&=-iz_2\Dz1-iz_4\Dz3+2w_3\Ew2,\\
X_{0,10}&=-z_1\Dz1+z_2\Dz2-z_3\Dz3+z_4\Dz4+2w_3\Ew3,\\
X_{0,11}&=z_1\Dz3+w_1\Ew4,\\
X_{0,12}&=iz_1\Dz3+w_2\Ew4,\\
X_{0,13}&=z_2\Dz3+2w_3\Ew4,\\
X_{0,14}&=z_3\Dz3+z_4\Dz4+w_4\Ew4.
\end{align*}

A basis of $g_1$:
\begin{align*}
X_{1,1}={}&(w_1+2iz_1z_2)\Dz3+2iz_2^2\Dz4+2iw_1z_2\Ew4,\\
X_{1,2}={}&(w_2-2z_1z_2)\Dz3+2z_2^2\Dz4+2iw_2z_2\Ew4,\\
X_{1,3}={}&(w_3+2iz_2^2)\Dz3+2iw_3z_2\Ew4,\\
X_{1,4}={}&2iz_1^2\Dz3+(w_1+2iz_1z_2)\Dz4+2iw_1z_1\Ew4,\\
X_{1,5}={}&-2z_1^2\Dz3+(w_2+2z_1z_2)\Dz4+2iw_2z_1\Ew4,\\
X_{1,6}={}&2iz_1z_2\Dz3+w_3\Dz4+2iw_3z_1\Ew4,\\[1mm]
X_{1,7}={}&(-W_{-}-4iz_1z_2)\Dz1+(-2w_3-4iz_2^2)\Dz2\\
&+(-2w_4-8iz_2z_3)\Dz3-8iz_2z_4\Dz4\\
&-2i(W_{+}z_2+2w_3z_1)\Ew1+(-2W_{+}z_2+4w_3z_1)\Ew2\\
&-4iw_3z_2\Ew3-2i(W_{+}z_4+2w_3z_3+2w_4z_2)\Ew4,\\[1mm]
X_{1,8}={}&(iw_1+2z_1z_2)\Dz3+2z_2^2\Dz4+2w_1z_2\Ew4,\\
X_{1,9}={}&(iw_2+2iz_1z_2)\Dz3-2iz_2^2\Dz4+2w_2z_2\Ew4,\\
X_{1,10}={}&(iw_3+2z_2^2)\Dz3+2w_3z_2\Ew4,\\[1mm]
X_{1,11}={}&(iW_{-}+4z_1z_2)\Dz1+(2iw_3+4z_2^2)\Dz2\\
&+(2iw_4+8z_2z_3)\Dz3+8z_2z_4\Dz4\\
&+2(W_{+}z_2+2w_3z_1)\Ew1+2i(-W_{+}z_2+2w_3z_1)\Ew2\\
&+4w_3z_2\Ew3+(2W_{+}z_4+4w_3z_3+4w_4z_2)\Ew4,\\[1mm]
X_{1,12}={}&2z_1^2\Dz3+(iw_1+2z_1z_2)\Dz4+2w_1z_1\Ew4,\\
X_{1,13}={}&2iz_1^2\Dz3+(iw_2-2iz_1z_2)\Dz4+2w_2z_1\Ew4,\\
X_{1,14}={}&2z_1z_2\Dz3+iw_3\Dz4+2w_3z_1\Ew4.
\end{align*}

A basis of $g_2$:
\begin{align*}
X_{2,1}={}&(iw_2z_2+w_3z_1)\Dz3-w_3z_2\Dz4,\\
X_{2,2}={}&iW_{+}z_1\Dz3+i(W_{-}z_2-4w_3z_1)\Dz4,\\
X_{2,3}={}&W_{+}z_1\Dz3-W_{-}z_2\Dz4,\\
X_{2,4}={}&i(-w_1z_2+w_3z_1)\Dz3+iw_3z_2\Dz4,\\
X_{2,5}={}&w_1z_1\Dz3+w_1z_2\Dz4+w_1^2\Ew4,\\
X_{2,6}={}&iw_1z_1\Dz3+(w_2z_2-2iw_3z_1)\Dz4+w_1w_2\Ew4,\\
X_{2,7}={}&(w_1+2iw_2)z_1\Dz3-w_1z_2\Dz4+w_2^2\Ew4,\\
X_{2,8}={}&(w_1z_2+w_3z_1)\Dz3+w_3z_2\Dz4+2w_1w_3\Ew4,\\
X_{2,9}={}&(-iw_1z_2+w_2z_2+2iw_3z_1)\Dz3+2w_2w_3\Ew4,\\
X_{2,10}={}&w_3z_2\Dz3+w_3^2\Ew4,\\[1mm]
X_{2,11}={}&W_{-}z_2\Dz1+2w_3z_2\Dz2\\
&+(-W_{+}z_4+2w_3z_3+2w_4z_2)\Dz3+4w_3z_4\Dz4\\
&+2w_1w_3\Ew1+2w_2w_3\Ew2+2w_3^2\Ew3+4w_3w_4\Ew4.
\end{align*}

A basis of $g_3$:
\begin{align*}
X_{3,1}={}&-iw_1^2+w_1w_2-2w_1z_1z_2+2iw_2z_1z_2+4w_3z_1^2\Dz3\\
&+2\,(iw_1w_3-w_1z_2^2+iw_2z_2^2+2w_3z_1z_2)\Dz4
+2\,w_1 T\Ew4,\\[1mm]
X_{3,2}={}&i(iw_1^2-w_1w_2-2w_1z_1z_2+2iw_2z_1z_2+4w_3z_1^2)\Dz3\\
&+2\,i(-iw_1w_3-w_1z_2^2+iw_2z_2^2+2w_3z_1z_2)\Dz4
+2\,iw_1T\Ew4,\\[1mm]
X_{3,3}={}&i(-w_1w_2-2w_1z_1z_2-iw_2^2+2iw_2z_1z_2+4w_3z_1^2)\Dz3\\
&-2\,i(-w_1z_2^2-w_2w_3+iw_2z_2^2+2w_3z_1z_2)\Dz4
+2\,w_2T\Ew4,\\[1mm]
X_{3,4}={}&-w_1w_2-2w_1z_1z_2+iw_2^2+2iw_2z_1z_2+4w_3z_1^2\Dz3\\
&+2\,(-w_1z_2^2+w_2w_3+iw_2z_2^2+2w_3z_1z_2)\Dz4
+2\,iw_2T\Ew4,\\[1mm]
X_{3,5}={}&-iw_1w_3-2w_1z_2^2+w_2w_3+2iw_2z_2^2+4w_3z_1z_2\Dz3\\
&+2\,iw_3^2\Dz4+2\,w_3T\Ew4,\\[1mm]
X_{3,6}={}&i(iw_1w_3-2w_1z_2^2-w_2w_3+2iw_2z_2^2+4w_3z_1z_2)\Dz3\\
&+2\,w_3^2\Dz4+2\,iw_3T\Ew4.
\end{align*}

The component $g_4$ is one-dimensional and is generated by the vector field
\begin{equation*}
X_{4,1}=i(w_1+iw_2)T\Dz3-2iw_3T\Dz4.
\end{equation*}
The component $g_5$ is zero. Since
$g_-=g_{-2}+g_{-1}$ is fundamental, all $g_j$ with $j\geq5$ vanish.

The resulting vector fields can be integrated explicitly without
difficulty, and the degrees of the resulting rational mappings do not
exceed three. The last assertion in part (b) follows from \cite{VB24}.
The proposition is proved.

\vspace{3ex}

For comparison, we give the profile of the Tanaka prolongation of Meylan's
example (see \cite{FM,GMs}), which served as a prototype for the
$(4,4)$-quadric described above:
\begin{equation*}
\begin{array}{c|rrrrrrrrr}
 j&-2&-1&0&1&2&3&4&5&6\\ \hline
\mathrm{dim}\, g_j &5&8&17&20&21&16&8&4&1
\end{array},
\end{equation*}
\[
\mathrm{dim} \, g = 5+8+17+20+21+16+8+4+1=100.
\]
In most known examples, the sequence of dimensions of the nonnegative
components decreases monotonically, and $g_0$ has maximal dimension. Here,
surprisingly, the strict maximum occurs at $g_2$. Also, the Levi
decomposition of $g$ has the form
$
 g\simeq \mathrm{su}(2,3) \ltimes
 \bigl(\RR \ltimes V\bigr),
$
where $V$ is an abelian real module of dimension $75$, while the additional
$\RR$ lies in $g_0$ and acts by scalars. Integrating the vector fields in
$\RR\ltimes V$ gives triangular polynomial transformations of degree at
most four. The transformations corresponding to the semisimple summand
$\mathrm{su}(2,3)$ are matrix fractional-linear mappings. An additional
argument shows that the overall degree bound does not exceed eight. The
result of \cite{VB24} applies to this local group. Hence its closure is an
algebraic subgroup of $Cr(\mathbf C^9)$.

\vspace{5ex}

{\bf 3. The Lattice of Types}

\vspace{3ex}

Let $\mathbf N=\{1,2,\dots\}$ and
$\mathbf L=\{(n,k)\in\mathbf N^2:k\leq n^2\}$, the lattice of types. A
lattice point $(n,k)\in\mathbf L$ is called exceptional if an exceptional
quadric of this type exists. Write $\mathbf L$ as the disjoint union
$A\sqcup B\sqcup C$, where\\
$A$ consists of the types currently proved to be nonexceptional (the
``nonexceptional zone'');\\
$B$ consists of the types for which an example is known, or its existence
has been proved (the ``exceptional zone'');\\
and $C$ is the ``gray zone'', consisting of the types for which neither an
example nor a nonexistence proof is currently known.

This partition will undoubtedly change over time as $C$ decreases and its
elements move either to $A$ or to $B$. Ideally, $C$ should eventually
disappear.

Our immediate goal is to describe these subsets of $\mathbf L$.

\vspace{3ex}

{\bf Proposition 2.} Let
\[Q=\{(z,w) \in \mathbf{C}^{n+k}: \mathrm{Im} \,w_j = <z,\bar{z}>_j,
\ j=1,\dots,k\}\]
be an exceptional $(n,k)$-quadric. Let $\zeta\in\mathbf C^\nu$, and let
$(\zeta,\bar\zeta)$ be a scalar nondegenerate Hermitian form on
$\mathbf C^\nu$. Then
\begin{eqnarray*}
 \tilde{Q}=\{(z,\zeta,w) \in \mathbf{C}^{n+\nu+k}: \mathrm{Im}\,w_j = <z,\bar{z}>_j, \; j=1,\dots , (k-1), \\
\mathrm{Im}\,w_k = <z,\bar{z}>_k+(\zeta,\bar{\zeta}) \} \qquad\qquad
\end{eqnarray*}
is an exceptional $(n+\nu,k)$-quadric.\\
{\it Proof.} The quadric $\tilde Q$ remains nondegenerate. Under the
natural identification, a nonzero vector field
$\mathcal X\in g_3(Q)$ becomes a vector field in $g_3(\tilde Q)$.

\vspace{3ex}

{\bf Corollary 3.} If $(n,k)\in B$ and $\nu\geq0$, then
$(n+\nu,k)\in B$.

\vspace{1ex}

Applying this construction to the exceptional $(4,4)$-quadric above, we
obtain, for every $n\geq4$, the quadric $Q_{n4}$ of type $(n,4)$
\begin{equation*}
\begin{aligned}
\mathrm{Im} w_1
&=z_1\overline z_2+z_2\overline z_1,\\
\mathrm{Im} w_2
&=iz_1\overline z_2-iz_2\overline z_1,\\
\mathrm{Im} w_3
&=|z_2|^2,\\
\mathrm{Im} w_4
&=z_1\overline z_4+z_4\overline z_1
  +z_2\overline z_3+z_3\overline z_2 + \sum_{j=5}^{n} |z_j|^2.
\end{aligned}
\end{equation*}
which is exceptional. Computing the full Tanaka prolongation of its symbol,
we obtain
\[
 g_{-2}\oplus  g_{-1}\oplus
 g_0\oplus g_1\oplus g_2
\oplus g_3\oplus g_4,
\]
where
\begin{equation*}
\begin{array}{c|r|r|r|r|r|r|r}
j&-2&-1&0&1&2&3&4\\
\hline
\mathrm{dim}\,  g_j&4&2n&n^2-4n+14&6 n-10&2n+3&6&1
\end{array} \; .
\end{equation*}

\vspace{3ex}

The cases $n=1$ and $n=2$ are well known and contain no exceptional types.
It was shown in \cite{VB21} that there are no exceptional types when
$n=3$ and $k\leq4$. The type $(3,5)$ is handled by similar elementary
arguments, while the types $(3,6),(3,7),(3,8),(3,9)$ are covered by
Proposition~7 below. Thus we obtain the following statement.

 \vspace{1ex}

{\bf Proposition 4.}
$A_0=\{(n,k)\in\mathbf L:\min(n,k)\leq3\}\subseteq A$.
  
   \vspace{1ex}

The direct product of quadrics of types $(n_1,k_1)$ and $(n_2,k_2)$ is a
quadric of type $(n_1+n_2,k_1+k_2)$. If both factors are nondegenerate,
then so is the product. If one factor is exceptional, then the product is
exceptional. Hence:

\vspace{1ex}

{\bf Proposition 5.} If $(n,k)\in B$ and
$(\nu,\kappa)\in\mathbf L$, then $(n+\nu,k+\kappa)\in B$.

 \vspace{1ex}   

Using the exceptional $(n,4)$-quadric described above and repeatedly
taking the nondegenerate $(1,1)$-quadric (the sphere) as a direct factor,
we obtain:

\vspace{1ex}   

{\bf Corollary 6.} $B_0=\{4\leq n,\ 4\leq k\leq n\}\subseteq B$.

\vspace{1ex}   

The exceptional zone $B$ in the lattice $\mathbf L$ can be enlarged
substantially by using the results of \cite{GM}.

\vspace{5ex}

The exceptional $(4,5)$-quadric there, together with Theorem~15, adds to
$B$ the entire ``superdiagonal'' $\{(n,n+1):n\geq4\}$. The quadrics in
Theorem~18 of that paper provide examples of still higher codimension.
They are notable because their graded Lie algebras have large, not
uniformly bounded, length, although this feature is not essential here.
The examples form two series. The first is parametrized by a positive
integer $m$ and contains exceptional quadrics of type $(2m,m^2+1)$. The
second is parametrized by a positive integer $r$ and contains exceptional
quadrics of type $(2r,r(r-1)/2+1)$. Applying the constructions above to
these quadrics gives
\[B_1=\{4 \leq n, \; 4 \leq k \leq (n-4)^2+5\} \subseteq B.\]
The only lattice point whose exceptionality follows from Theorem~18 but is
not contained in $B_1$ is then $(6,10)$.

\vspace{3ex}

We now turn to the nonexceptional zone $A$. What does $A$ contain in
addition to $A_0$?

\vspace{3ex}
{\bf Proposition 7.} If $(n,k)\in\mathbf L$ and $k\geq n^2-n$, then
$(n,k)\in A$.\\
{\it Proof.} Consider the complex vector space
\[
 <H>=\mathrm{span}_{\mathbb C}\{H_1,\ldots,H_k\}
 \subset\operatorname{Mat}_n(\mathbb C),\qquad
 \dim_{\mathbb C}<H>=k,
\]
and put
\[
 \mathcal N=<H>^\perp
 =\{N:\operatorname{tr}(H_jN)=0,\ j=1,\ldots,k\}.
\]
Then $\dim_{\mathbb C}\mathcal N=n^2-k$. For $a\in\mathbf C^n$, set
$a^H=\{x\in\mathbf C^n:<x,a>=0\}$. The condition $x\in a^H$ means
\[
 \langle H_jx,a\rangle=a^*H_jx=0,
 \qquad j=1,\ldots,k,
\]
while $a^*H_jx=\operatorname{tr}(H_jxa^*)$. Consequently,
\[
 x\in a^H\quad\Longleftrightarrow\quad xa^*\in\mathcal N.
\]
For $x\ne0$, the matrix $R_a=xa^*$ has rank one and is a decomposable
tensor $x\otimes a^*$. Thus $a^H\ne0$ exactly when $\mathcal N$ contains a
nonzero matrix $x\otimes a^*$ with the prescribed right factor $a^*$.
The space of all such matrices has dimension $n$. Hence, if
$k>n^2-n$, we conclude that $a^H=0$. By Proposition~6.4 of \cite{T}, the
collection $H=(H_1,\dots,H_k)$ is then $T$-nondegenerate, and Theorem~5.2
of that paper implies that both the form $H$ and the quadric $Q(H)$ are
nonexceptional.

It remains to consider $k=n^2-n$. Then
$\dim\mathcal N=n$. If $\mathcal N$ consisted entirely of matrices of
rank one, this linear space would have either a common one-dimensional
image or a common kernel of codimension one. For $n>1$, neither type of
subspace can be invariant under Hermitian conjugation, a contradiction.
It follows that there is an $a$ for which\footnote{The Russian source
prints $k$ in the following formula. Since the displayed matrix has $n$
rows, $n$ is the mathematically possible full-rank condition.}
\[
 \operatorname{rank}_{\mathbf C}(H_1a,\dots,H_ka)=n.
\]
Proposition~6.3 and Theorem~5.2 of Tumanov \cite{T} now give $g_3=0$.
The proposition is proved.

\vspace{3ex}

{\bf Theorem 8.}
The set $A$ of nonexceptional types contains
\[
 \{(n,k)\in\mathbf L:n\leq3\}\cup
 \{(n,k)\in\mathbf L:k\leq3\}\cup
 \{(n,k)\in\mathbf L:k\geq n^2-n\}.
\]
The set $B$ of exceptional types contains
\[
 \{(n,k)\in\mathbf L:4\leq n,\ 4\leq k\leq(n-4)^2+5\}
 \cup\{(6,10)\}.
\]

\vspace{3ex}

All other types currently belong to the gray zone $C$. The smallest type,
in the sense of minimizing $n+k$, occurring in $C$ is $(4,6)$. Apart from
the point $(6,10)$, the set $C$ consists of the lattice points between two
parabolas:
\[
C=\{(n,k)\in\mathbf L:4\leq n,\ (n-4)^2+6\leq k\leq n^2-n-1\}
\setminus\{(6,10)\}.
\]

Let $a(n)$, $b(n)$, and $c(n)$ be the numbers of types with fixed $n$ in
$A$, $B$, and $C$, respectively. If $n\geq7$, then
$a(n)=n+4$, $b(n)=(n-4)^2+1$, and $c(n)=7n-22$. Hence the proportion of
$A\cup C$ tends to zero as $n$ grows, and almost all types belong to $B$.
Thus exceptional types occur very frequently. What is rare, or
exceptional, is not a type admitting an exceptional quadric, but the
exceptional quadric itself.
\begin{figure}[htbp]
\centering
\includegraphics[width=0.92\linewidth]
  {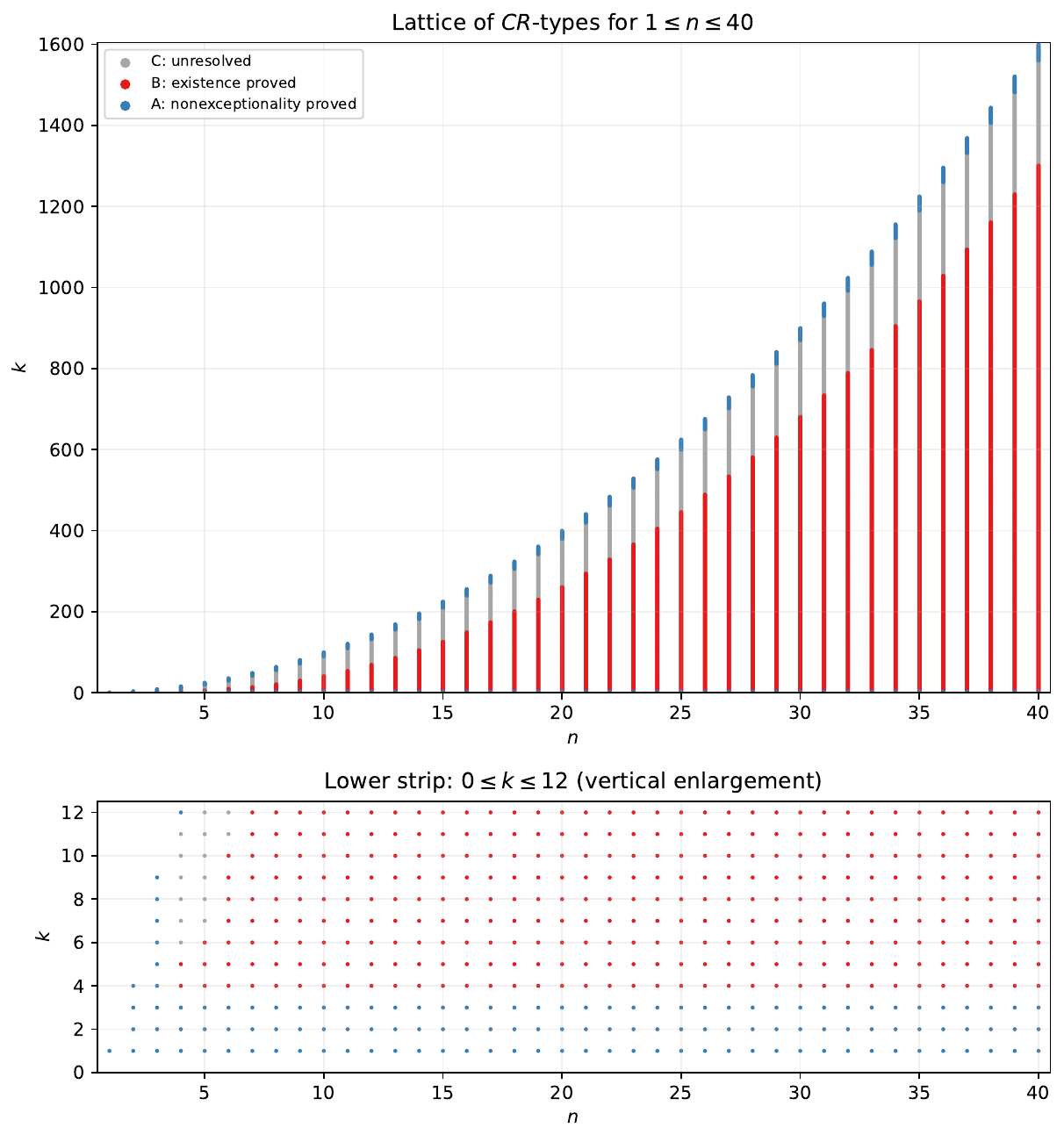}
\end{figure}

\vspace{3ex}

For the $(4,4)$-quadric constructed above, $l_+$, the length of $g_+$, is
equal to $4$. The general estimate in \cite{VB21} implies that
$l_+\leq8$ for every quadric of codimension four. This leads to the
following question.

\vspace{2ex}

{\bf Question 9.} What is $\max l_+(Q)$ over all exceptional
$(4,4)$-quadrics $Q$?

\vspace{2ex}
The smallest dimension of an ambient space
$\mathbf C^N=\mathbf C^{n+k}$ containing an unresolved type, that is, a
type in $C$, is $N=10$. There is one such type, namely $(4,6)$.
  
\vspace{2ex}

{\bf Question 10.} Do exceptional $(4,6)$-quadrics exist?

\vspace{2ex}

Maple and ChatGPT were used in the preparation of this paper.

\end{document}